\documentclass[11pt]{amsart}

\usepackage[mathscr]{eucal}
\usepackage{amsfonts}
\usepackage{amsmath}
\usepackage{amsthm}
\usepackage{amssymb}
\usepackage{latexsym}
\usepackage{pstricks,pst-node}

\theoremstyle{plain}

\newtheorem{Thm}[equation]{Theorem}

\renewcommand{\i}{${\sf {(i)}}\;$}
\newcommand{\ii}{${\sf {(ii)}}\;$}
\newcommand{\iii}{${\sf {(iii)}}\;$}
\newcommand{\iv}{${\sf {(iv)}}\;$}
\newcommand{\vv}{${\sf {(v)}}\;$}
\renewcommand{\a}{\alpha}
\newcommand{\G}{\mathfrak{G}}
\newcommand{\Id}{\mathrm{Id}}
\renewcommand{\o}{\mathsf{or}}
\newcommand{\<}{\langle}
\renewcommand{\>}{\rangle}
\renewcommand{\v}{^\vee}

\begin{document}

\author{Wee Liang Gan}
\address{Department of Mathematics, Massachusetts Institute
of Technology, Cambridge, MA 02139, U.S.A.}
\email{wlgan@math.mit.edu}  

\title{On a theorem of Conant-Vogtmann}

\begin{abstract}
We prove that the graph complex is a strong homotopy Lie 
(super) bialgebra. 
\end{abstract}

\maketitle

\section{Introduction and definitions}

The graph complex was introduced by Kontsevich in \cite{K1} and
\cite{K2}. In \cite{CV1}, Conant and Vogtmann constructed a new 
differential, a Lie bracket and a Lie cobracket on the 
graph complex. They proved that on the subspace of one-particle
irreducible graphs, the Lie bracket and the Lie cobracket
give a Lie bialgebra structure. In this paper, we prove that
the whole graph complex is a \emph{strong homotopy Lie (super) 
bialgebra} (cf. \cite[\S 5]{M}, \cite[\S 3.2]{G}).

We shall work in the category of super vector spaces 
$V=V_+\oplus V_-$ over a field of characteristic $0$. 
Denote by $\Pi$ the parity-change functor $(\Pi V)_\pm = V_\mp$.

By a graph, we mean a finite $1$-dimensional CW-complex
(not necessarily connected) such that the valency of each 
vertex is at least $3$. An orientation $\o$ of a graph $X$ is 
an ordering of its vertices and a choice of direction on each
of its edges (cf. \cite[Definition 2]{CV2}). 
Switching the order of two vertices or reversing
the direction of an edge changes the orientation from $\o$ to
$-\o$. We shall identify an ordering of the vertices with a
labeling of the vertices by $1, 2, \ldots$.

Let $\G$ be the vector space generated by graphs $(X, \o)$,
modulo the equivalence relation $(X, -\o)= -(X,\o)$.
Observe that if $X$ has an edge-loop, then it is equal to $0$
in $\G$.
Let $\G_+$ (resp. $\G_-$) be the subspace of $\G$ spanned by
graphs with an even (resp. odd) number of vertices.
Clearly, $\G=\G_+\oplus\G_-$.
If $X,Y\in\G$, then define their product $X\cdot Y$ to be the
disjoint union of $X$ and $Y$. The orientation
of $X\cdot Y$ is obtained by adding the number of vertices of
$X$ to the labels of the vertices of $Y$. The vector space
$\G$ is a commutative super algebra by extending this product
bilinearly to linear combinations of graphs.

Given an edge $e$ of a graph $X\in\G$ such that $e$ is not 
an edge-loop, we denote by $X_e$ the
graph obtained from $X$ by contracting the edge $e$. 
The induced orientation on $X_e$ is defined as follows: 
if the source vertex of $e$ is labeled $1$ and
the target vertex of $e$ is labeled $2$, then the vertex
of $X_e$ formed by the contracted edge $e$ is labelled $1$,
and the labeling of other vertices are reduced by $1$.
If $e$ is an edge-loop, then we set $X_e$ to be $0$.

Given a half-edge $h$ of a graph $X\in\G$, we denote by $v(h)$ 
the vertex of $X$ attached to $h$, and by $\bar{h}$ the
half-edge such that $e(h):=h\cup\bar{h}$ is an edge of $X$.
Let $h_1, h_2$ be two half-edges such that they belong to two 
distinct edges of $X$, and choose a representative for the 
orientation on $X$ so that $v(h_1)$ (resp. $v(h_2)$)
is the source vertex of $e(h_1)$ (resp. $e(h_2)$).
Assume that $e(h_1)$ and $e(h_2)$ are not edge-loops.
Denote by $X\< h_1,h_2\>$ the graph obtained from $X$ by
adding in an edge $e_1$ directed from $v(h_1)$ to $v(\bar{h}_2)$, 
adding in an edge $e_2$ directed from $v(h_2)$ to $v(\bar{h}_1)$,
and deleting away the two edges $e(h_1), e(h_2)$.
Denote by $h'_i$ ($i=1,2$) the half-edge of $X\< h_1,h_2\>$
such that $e_i=e(h'_i)$ and $e_i$ is directed from
$h'_i$ to $\bar{h}'_i$; see following picture:
$$
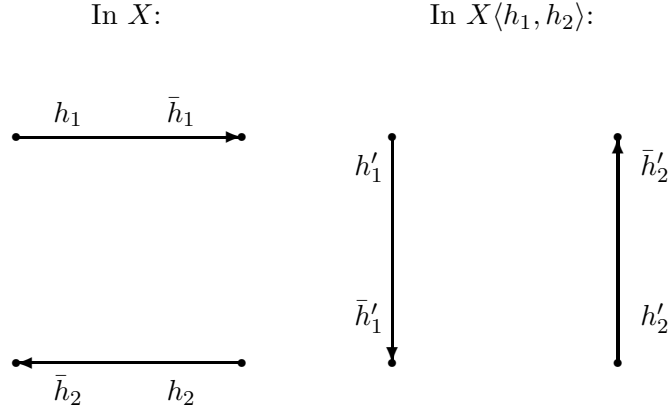
\setlength{\unitlength}{1mm}
\begin{picture}(100,60)
\put(10,40){\circle*{1}}
\put(40,40){\circle*{1}}
\put(10,10){\circle*{1}}
\put(40,10){\circle*{1}}
\put(60,40){\circle*{1}}
\put(90,40){\circle*{1}}
\put(60,10){\circle*{1}}
\put(90,10){\circle*{1}}
\thicklines
\put(10,40){\vector(1,0){30}}
\put(40,10){\vector(-1,0){30}}
\put(60,40){\vector(0,-1){30}}
\put(90,10){\vector(0,1){30}}
\put(20,55){$\mbox{In $X$:}$}
\put(65,55){$\mbox{In $X\< h_1,h_2\>$:}$}
\put(15,42){$h_1$}
\put(30,42){$\bar h_1$}
\put(15,5){$\bar h_2$}
\put(30,5){$h_2$}
\put(55,35){$h'_1$}
\put(55,15){$\bar h'_1$}
\put(93,35){$\bar h'_2$}
\put(93,15){$h'_2$}
\end{picture}
$$
Finally, denote $X_{h_1, h_2}:= X\< h_1,h_2\>_{e_1}$.
If $v$ is a vertex of $X\< h_1,h_2\>$ attached to $e_1$, 
we will also denote by $v$ the vertex of $X_{h_1, h_2}$ 
contracted from the edge $e_1$. (Thus, the vertex of $X_{h_1, h_2}$
contracted from $e_1$ may be denoted by more than one name.)
Observe that $X_{\bar{h}_2, \bar{h}_1} = X_{h_1, h_2}$.
If $e(h_1)$ or $e(h_2)$ is an edge-loop, then we set
$X_{h_1, h_2}$ to be $0$. (Our notation differs from the one
in \cite{CV1}, where they write $X_{h_1\bar{h_2}}$ for our
$X_{h_1, h_2}$.)

For positive integers $m,n\leq 2$, we define an odd linear 
map $\a_{m,n}:\G^{\otimes n} \to \G^{\otimes m}$ as follows:
\begin{equation} \label{eqn1}
\a_{m,n}(X_1 \otimes \cdots \otimes X_n) :=
\sum_{h_1, h_2} \sum_{Y_1 \otimes \cdots \otimes
Y_m} Y_1 \otimes \cdots \otimes Y_m.
\end{equation}
In (\ref{eqn1}), the first summation is taken over all
ordered pairs of half-edges $h_1, h_2$ 
of $X_{1} \cdots X_{n}$ such that each graph $X_{i}$
contains $h_j$ for some $j$, and $h_1, h_2$ belong to
distinct edges. 
The second summation is taken over all $Y_{1} \otimes \cdots
\otimes Y_{m} \in\G^{\otimes m}$ such that
\begin{equation} \label{eqn02} 
Y_{1} \cdots Y_{m} = (X_{1} \cdots X_{n})_{h_1, h_2},
\end{equation} 
and each $Y_{i}$ is a graph containing $v(h'_j)$ for some $j$. 
If $m$ and $n$ are positive integers such that
$m>2$ or $n>2$, then we let
$\a_{m,n} : \G^{\otimes n} \to \G^{\otimes m}$
be the zero map.
The maps $\a_{1,1}$, $\a_{1,2}$, and $\a_{2,1}$ are, respectively,
the differential, Lie bracket, and Lie cobracket constructed 
in \cite{CV1}, but we will not use these facts in the proof of
our theorem below. The map $\a_{2,2}$ is new.

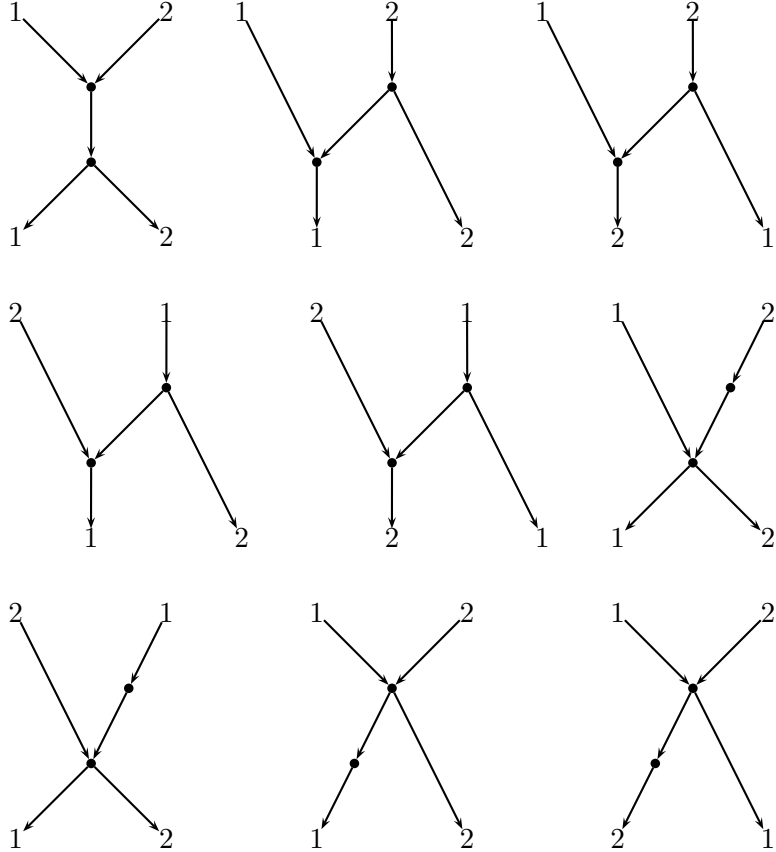
\begin{figure} \label{ff}

\begin{pspicture}(0,0)(10,11)

\rput(0,11){\rnode{A}{1}}
\rput(2,11){\rnode{B}{2}}
\dotnode(1,10){C}
\dotnode(1,9){D}
\rput(0,8){\rnode{E}{1}}
\rput(2,8){\rnode{F}{2}}
\ncline{->}{A}{C} 
\ncline{->}{B}{C} 
\ncline{->}{C}{D}
\ncline{->}{D}{E}
\ncline{->}{D}{F}

\rput(3,11){\rnode{A}{1}}
\rput(5,11){\rnode{B}{2}}
\dotnode(5,10){C}
\dotnode(4,9){D}
\rput(4,8){\rnode{E}{1}}
\rput(6,8){\rnode{F}{2}}
\ncline{->}{A}{D}
\ncline{->}{B}{C}
\ncline{->}{C}{D}
\ncline{->}{D}{E}
\ncline{->}{C}{F}

\rput(7,11){\rnode{A}{1}}
\rput(9,11){\rnode{B}{2}}
\dotnode(9,10){C}
\dotnode(8,9){D}
\rput(8,8){\rnode{E}{2}}
\rput(10,8){\rnode{F}{1}}
\ncline{->}{A}{D}
\ncline{->}{B}{C}
\ncline{->}{C}{D}
\ncline{->}{D}{E}
\ncline{->}{C}{F}

\rput(0,7){\rnode{A}{2}}
\rput(2,7){\rnode{B}{1}}
\dotnode(2,6){C}
\dotnode(1,5){D}
\rput(1,4){\rnode{E}{1}}
\rput(3,4){\rnode{F}{2}}
\ncline{->}{A}{D}
\ncline{->}{B}{C}
\ncline{->}{C}{D}
\ncline{->}{D}{E}
\ncline{->}{C}{F}

\rput(4,7){\rnode{A}{2}}
\rput(6,7){\rnode{B}{1}}
\dotnode(6,6){C}
\dotnode(5,5){D}
\rput(5,4){\rnode{E}{2}}
\rput(7,4){\rnode{F}{1}}
\ncline{->}{A}{D}
\ncline{->}{B}{C}
\ncline{->}{C}{D}
\ncline{->}{D}{E}
\ncline{->}{C}{F}

\rput(8,7){\rnode{A}{1}}
\rput(10,7){\rnode{B}{2}}
\dotnode(9.5,6){C}
\dotnode(9,5){D}
\rput(8,4){\rnode{E}{1}}
\rput(10,4){\rnode{F}{2}}
\ncline{->}{A}{D}
\ncline{->}{B}{C}
\ncline{->}{C}{D}
\ncline{->}{D}{E}
\ncline{->}{D}{F}

\rput(0,3){\rnode{A}{2}}
\rput(2,3){\rnode{B}{1}}
\dotnode(1.5,2){C}
\dotnode(1,1){D}
\rput(0,0){\rnode{E}{1}}
\rput(2,0){\rnode{F}{2}}
\ncline{->}{A}{D}
\ncline{->}{B}{C}
\ncline{->}{C}{D}
\ncline{->}{D}{E}
\ncline{->}{D}{F}

\rput(4,3){\rnode{A}{1}}
\rput(6,3){\rnode{B}{2}}
\dotnode(5,2){C}
\dotnode(4.5,1){D}
\rput(4,0){\rnode{E}{1}}
\rput(6,0){\rnode{F}{2}}
\ncline{->}{A}{C}
\ncline{->}{B}{C}
\ncline{->}{C}{D}
\ncline{->}{D}{E}
\ncline{->}{C}{F}

\rput(8,3){\rnode{A}{1}}
\rput(10,3){\rnode{B}{2}}
\dotnode(9,2){C}
\dotnode(8.5,1){D}
\rput(8,0){\rnode{E}{2}}
\rput(10,0){\rnode{F}{1}}
\ncline{->}{A}{C}
\ncline{->}{B}{C}
\ncline{->}{C}{D}
\ncline{->}{D}{E}
\ncline{->}{C}{F}

\end{pspicture}

\caption{Elements of $T(2,2)$.}
\end{figure}

By a corolla, we mean a vertex with directed half-edges attached
to it such that there is at least one incoming half-edge and
at least one outgoing half-edge.
Let $m,n$ be positive integers. We shall denote by $T(m,n)$
the set consisting of all ``flowcharts'' $T$ described as follows:
$T$ is obtained from two corollas, called $s(T)$ and $t(T)$,
by joining an outgoing half-edge of $s(T)$ to an incoming 
half-edge of $t(T)$, such that the resulting flowchart has
$n$ inputs and $m$ outputs, and moreover there is a labeling
of the inputs by $1, \ldots, n$ and a labeling of the outputs 
by $1, \ldots, m$.
For example, the elements of $T(2,2)$ are listed in Figure
\ref{ff}.
If $v$ is a corolla, we write $i(v)$ (resp. $o(v)$) for the 
number of incoming (resp. outgoing) half-edges at $v$.
For each flowchart $T\in T(m,n)$, we write 
$$ \a_{o(t(T)),i(t(T))} \circ_T \a_{o(s(T)),i(s(T))} :
\G^{\otimes n} \to \G^{\otimes m}$$
for the composition of $\a_{o(t(T)),i(t(T))}$ and 
$\a_{o(s(T)),i(s(T))}$ according to $T$.

\begin{Thm} \label{main}
For each ordered pair of positive integers $m,n$, 
the map $\a_{m,n}$ is commutative and cocommutative, and
the following identity holds:
\begin{equation} \label{eqn11}
\sum_{T\in T(m,n)} \a_{o(t(T)),i(t(T))} \circ_T 
\a_{o(s(T)),i(s(T))} = 0.
\end{equation}
In other words, $\Pi\G$ is a strong homotopy Lie (super) bialgebra.
\end{Thm}

We will prove Theorem \ref{main} in section \ref{s2}. 
In the appendix at the end of the paper, we 
give the definition of strong homotopy Lie bialgebra.

Note that the left hand side of (\ref{eqn11}) is zero if 
$m$ or $n$ is greater than $3$ because $\a_{i,j}=0$ if 
$i$ or $j$ is greater than $2$.
Taking $(m,n)=(1,1)$, $(1,2)$, $(1,3)$, $(2,1)$, $(3,1)$
in (\ref{eqn11}), we obtain the results of \cite{CV1} that
$\a_{1,1}$ is a differential, 
$\a_{1,1}$ is a derivation with respect to $\a_{1,2}$,
$\a_{1,2}$ satisfies the Jacobi identity, 
$\a_{1,1}$ is a coderivation with respect to $\a_{2,1}$, and 
$\a_{2,1}$ satisfies the coJacobi identity.
Taking $(m,n)=(2,2)$, $(2,3)$, $(3,2)$ and $(3,3)$, 
we obtain new identities.
In particular, taking $(m,n)=(2,2)$, we deduce
\cite[Theorem 1]{CV1} which states that $\Pi\G^{1\textrm{PI}}$
is a Lie super bialgebra, where $\G^{1\textrm{PI}} \subset\G$ 
is the subspace spanned by one-particle irreducible 
graphs (i.e. connected graphs which remain connected after 
the removal of any edge).
Indeed, it is plain that $\a_{2,2}(X_1\otimes X_2)=0$ if
$X_1$ or $X_2$ is in $\G^{1\textrm{PI}}$.

We remark that in \cite{C}, Conant constructed a 
family of strong homotopy Lie algebra structures
and strong homotopy Lie coalgebra structures
on graph complexes associated to cyclic operads. 
We do not know if these structures also form strong 
homotopy Lie bialgebras.

\section{Proof of Theorem \ref{main}} \label{s2}

The commutativity and cocommutativity of $\a_{m,n}$ are clear
from (\ref{eqn02}). We have to prove that (\ref{eqn11}) holds.
To this end, we shall pair the terms which appear in the
left hand side of (\ref{eqn11}) such that the terms in
a pair are negative of each other.
Roughly speaking, the terms which come from two pairs
of half-edges taken in different order have opposite signs.
However, this is not always true, as we now illustrate:
$$ 
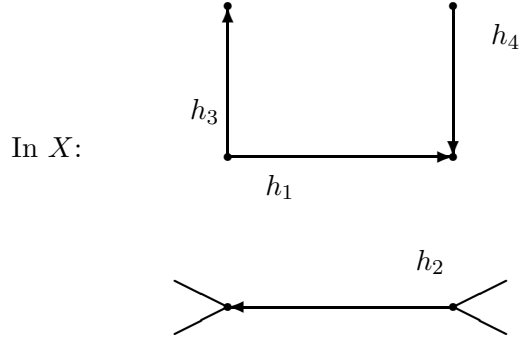
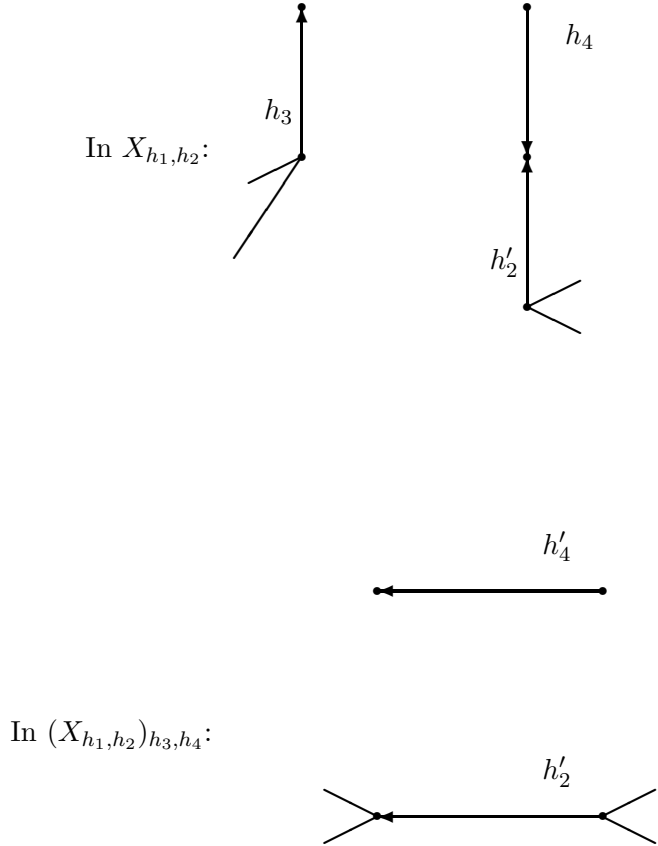
\setlength{\unitlength}{1mm}
\begin{picture}(80,50)
\put(30,50){\circle*{1}}
\put(30,30){\circle*{1}}
\put(30,10){\circle*{1}}
\put(60,50){\circle*{1}}
\put(60,30){\circle*{1}}
\put(60,10){\circle*{1}}
\thicklines
\put(30,30){\vector(0,1){20}}
\put(30,30){\vector(1,0){30}}
\put(60,50){\vector(0,-1){20}}
\put(60,10){\vector(-1,0){30}}
\put(35,25){$h_1$}
\put(55,15){$h_2$}
\put(25,35){$h_3$}
\put(65,45){$h_4$}
\put(0,30){$\mbox{ In $X$: }$}
\put(30,10){\line(-2,1){7}}
\put(30,10){\line(-2,-1){7}}
\put(60,10){\line(2,1){7}}
\put(60,10){\line(2,-1){7}}
\end{picture}
$$

$$
\setlength{\unitlength}{1mm}
\begin{picture}(80,50)
\put(30,50){\circle*{1}}
\put(30,30){\circle*{1}}
\put(60,10){\circle*{1}}
\put(60,50){\circle*{1}}
\put(60,30){\circle*{1}}
\thicklines
\put(30,30){\vector(0,1){20}}
\put(60,50){\vector(0,-1){20}}
\put(60,10){\vector(0,1){20}}
\put(55,15){$h'_2$}
\put(25,35){$h_3$}
\put(65,45){$h_4$}
\put(0,30){$\mbox{ In $X_{h_1,h_2}$: }$}
\put(30,30){\line(-2,-1){7}}
\put(30,30){\line(-2,-3){9}}
\put(60,10){\line(2,1){7}}
\put(60,10){\line(2,-1){7}}
\end{picture}
$$

$$
\setlength{\unitlength}{1mm}
\begin{picture}(100,70)
\put(50,50){\circle*{1}}
\put(50,20){\circle*{1}}
\put(80,20){\circle*{1}}
\put(80,50){\circle*{1}}
\thicklines
\put(80,50){\vector(-1,0){30}}
\put(80,20){\vector(-1,0){30}}
\put(72,25){$h'_2$}
\put(72,55){$h'_4$}
\put(0,30){$\mbox{ In $(X_{h_1,h_2})_{h_3,h_4}$: }$}
\put(50,20){\line(-2,1){7}}
\put(50,20){\line(-2,-1){7}}
\put(80,20){\line(2,1){7}}
\put(80,20){\line(2,-1){7}}
\end{picture}
$$
However, in this example, the edge $e(h_1)$ becomes an edge-loop in
$X_{h_3,h_4}$ and hence $X_{h_3,h_4}$ is equal to zero.
In below, we shall pair the term $(X_{h_1,h_2})_{h_3,h_4}$ with
$(X_{\bar h_1,h_2})_{ h_3, h_4}$, which is equal to
$-(X_{h_1,h_2})_{h_3,h_4}$ (the minus sign is due to
orientation).

To make precise how the terms on the left hand side of 
(\ref{eqn11}) cancel, 
we shall define a set $F$ and an involution
$\mu : F \to F$. The terms in the left hand side of
(\ref{eqn11}) will be grouped according to $F$, and the terms
corresponding to $f\in F$ will cancel with the terms
corresponding to $\mu(f)\in F$.

Fix any $m,n\leq 3$ and graphs $X_1, \ldots, X_n\in \G$. 
We may assume that there is no edge-loop in $X_1\cdots X_n$.
Define $F$ to be the set consisting of all data
$$ f:= (h_1, h_2, h_3, h_4, U_1, \ldots, U_m) $$
where
\begin{itemize}
\item
$h_1, h_2$ are half-edges of $X_1\cdots X_n$ which belong to 
distinct edges;
\item
$e(h'_1), e(h'_2)$ are not edge-loops in $(X_1\cdots X_n)\<h_1,h_2\>$;
\item
$h_3, h_4$ are half-edges of $(X_1\cdots X_n)_{h_1,h_2}$
which belong to distinct edges;
\item
$e(h_3), e(h_4)$ are not edge-loops in $(X_1\cdots X_n)_{h_1,h_2}$, 
and $e(h'_3), e(h'_4)$ are not edge-loops in
$(X_1\cdots X_n)_{h_1,h_2}\< h_3,h_4\>$;
\item
defining $S_i$ ($i=1,\ldots, n$)
to be the set of all $j\in \{1,2,3,4\}$
such that either 
\begin{itemize}
\item $X_i$ contains $h_j$, or
\item $X_i$ contains $h_2$, and $h_j=h'_2$, or
\item $X_i$ contains $h_1$, and $h_j=\overline{h'_2}$, 
\end{itemize}
we have:
each $S_i$ is nonempty and 
$S_1\cup\ldots\cup S_n = \{1,2,3,4\}$ is a disjoint union;
\item
at least one $S_i$ which contains $1$ or $2$ also contains
$3$ or $4$;
\item
each $U_i$ is a nonempty subset of $\{1,2,3,4\}$ and
$U_1\cup\ldots\cup U_m = \{1,2,3,4\}$ is a disjoint union;
\item
at least one $U_i$ which contains $3$ or $4$ also contains
$1$ or $2$;
\item
if $v(h'_j)$, $v(h'_k)$ are in the same connected component of
$$((X_1\cdots X_n)_{h_1,h_2})_{h_3,h_4}$$ and $U_i$ contains $j$,
then $U_i$ contains $k$.
\end{itemize}
Given $f\in F$ as above, there is a unique flowchart 
$T_f\in T(m,n)$ such that the input of $T_f$ labelled $i$
goes into $s(T_f)$ if and only if $S_i$ contains $1$ or $2$,
and the output of $T_f$ labelled $i$ goes out from $t(T_f)$ 
if and only if $U_i$ contains $3$ or $4$.

Let $m,n \leq 3$.
Consider the left hand side of (\ref{eqn11}) applied to
$X_1\otimes\cdots\otimes X_n$. 
For each $(h_1,h_2,h_3,h_4, U_1, \ldots, U_m)\in F$,
there is a corresponding term
\begin{equation} \label{eqn5}
\sum_{Y_1\otimes\cdots\otimes Y_m} Y_1\otimes\cdots\otimes Y_m
\end{equation}
where the summation is taken over all 
$Y_1\otimes\cdots\otimes Y_m \in\G^{\otimes m}$ such that
$Y_1\cdots Y_m = ((X_1\cdots X_n)_{h_1,h_2})_{h_3,h_4}$
and $Y_i$ is a graph containing $v(h'_j)$ if $U_i$ contains $j$.

We define an involution $\mu:F\to F$ by
$$ \mu(h_1, h_2, h_3, h_4, U_1, \ldots, U_m)
:= (h_1\v, h_2\v, h_3\v, h_4\v, U_1\v, \ldots, U_m\v),$$
$$ U_i\v := \{j \mid \xi(j)\in U_i\},\quad
\xi: \{1,2,3,4\} \to \{1,2,3,4\},$$
where $h_1\v, h_2\v, h_3\v, h_4\v$ and 
$U_1\v, \ldots, U_m\v$ will be defined
accordingly in each of the following cases:

Case \i, $h_3,h_4 \nsubseteq e(h_2)$:
\begin{itemize}
\item 
if $e(h_1)$ connects $v(h_3)$ with $v(\bar{h}_4)$, then
$$(h_1\v, h_2\v, h_3\v, h_4\v) = (\bar{h}_1, h_2, h_3, h_4),
\quad \xi=\Id;$$
\item
if $e(h_2)$ connects $v(h_3)$ with $v(\bar{h}_4)$, then
$$(h_1\v, h_2\v, h_3\v, h_4\v) = (h_1, \bar{h}_2, h_3, h_4),
\quad \xi=\Id;$$
\item
otherwise, let
$$(h_1\v, h_2\v, h_3\v, h_4\v) = (h_3,h_4,h_1,h_2),$$
$$\xi(1)=3,\quad \xi(2)=4,\quad \xi(3)=1,\quad \xi(4)=2.$$
\end{itemize}

Case \ii, $h_3=h'_2$:
\begin{itemize}
\item
if $v(h_1)=v(h_2)$ and $v(\bar{h}_1)=v(\bar{h}_4)$, then
$$(h_1\v, h_2\v, h_3\v) = (h_1, h_4, \bar{h}_2),
\quad h_4\v = {h_2\v}', \quad \xi=\Id;$$
\item
if $v(\bar{h}_2)=v(h_4)$ and $v(h_1)=v(h_2)$, then
$$(h_1\v, h_2\v, h_3\v) = (\bar{h}_4, h_2, h_1),
\quad h_4\v = {h_2\v}', \quad \xi=\Id;$$
\item
if $v(\bar{h}_2)=v(h_4)$ and $v(\bar{h}_1)=v(\bar{h}_4)$, then
$$(h_1\v, h_2\v, h_3\v) = (\bar{h}_1, \bar{h}_2, h_4), 
\quad h_4\v = {h_2\v}', \quad \xi=\Id;$$
\item
if $v(\bar{h}_2)=v(h_4)$, $v(h_1)\neq v(h_2)$ and
$v(\bar{h}_1)\neq v(\bar{h}_4)$, then
$$(h_1\v, h_2\v, h_3\v) = (h_4, \bar{h}_1, h_2),
\quad h_4\v = {h_2\v}', $$
$$ \xi(1)=1,\quad \xi(2)=1,\quad \xi(3)=3,\quad \xi(4)=4;$$
\item
otherwise, let
$$ (h_1\v, h_2\v, h_3\v) = (h_2, h_4, h_1),
\quad h_4\v = {h_2\v}',$$
$$\xi(1)=2,\quad \xi(2)=4,\quad \xi(3)=1,\quad \xi(4)=4.$$
\end{itemize}

Case \iii, $h_4=h'_2$:
\begin{itemize}
\item
if $v(h_2)=v(h_3)$ and $v(\bar{h}_2)=v(\bar{h}_1)$, then
$$(h_1\v, h_2\v, h_4\v) = (h_3, h_2, \bar{h}_1),
\quad h_3\v = {h_2\v}' , \quad \xi=\Id;$$
\item
if $v(\bar{h}_3)=v(h_1)$ and $v(\bar{h}_2)=v(\bar{h}_1)$, then
$$(h_1\v, h_2\v, h_4\v) = (h_1, \bar{h}_3, h_2),
\quad h_3\v = {h_2\v}' , \quad \xi=\Id;$$
\item
if $v(\bar{h}_3)=v(h_1)$ and $v(h_2)=v(h_3)$, then
$$(h_1\v, h_2\v, h_4\v) = (\bar{h}_1, \bar{h}_2, h_3),
\quad h_3\v = {h_2\v}' , \quad \xi=\Id;$$
\item
if $v(\bar{h}_3)=v(h_1)$, $v(\bar{h}_2)\neq v(\bar{h}_1)$
and $v(h_2)\neq v(h_3)$, then 
$$(h_1\v, h_2\v, h_4\v) = (\bar{h}_2, h_3, h_1),
\quad h_3\v = {h_2\v}' ,$$
$$ \xi(1)=1,\quad \xi(2)=3,\quad \xi(3)=3, \quad \xi(4)=4;$$
\item
otherwise, let
$$(h_1\v, h_2\v, h_4\v) = (h_3, h_1, h_2),
\quad h_3\v = {h_2\v}' ,$$
$$\xi(1)=3,\quad \xi(2)=1,\quad \xi(3)=1,\quad \xi(4)=2.$$
\end{itemize}

Case \iv, $h_3=\overline{h'_2}$: 
$$ (h_1\v, h_2\v, \overline{h_4\v}, \overline{h_3\v}, 
U_1\v, \ldots, U_m\v) = 
\mu(h_1, h_2, \bar{h}_4, \bar{h}_3, U_1, \ldots, U_m)$$
where the right hand side is defined by case \iii.

Case \vv, $h_4=\overline{h'_2}$: 
$$ (h_1\v, h_2\v, \overline{h_4\v}, \overline{h_3\v},
U_1\v, \ldots, U_m\v) =
\mu(h_1, h_2, \bar{h}_4, \bar{h}_3, U_1, \ldots, U_m)$$
where the right hand side is defined by case \ii.

In the above, ${h_2\v}'$ is always taken to be in
$(X_1\cdots X_n)_{h_1\v, h_2\v}$.
It is a straightforward check that $\mu$ defines a pairing between
elements of $F$, and moreover,
$$((X_1\cdots X_n)_{h_1, h_2})_{h_3,h_4} 
= -((X_1\cdots X_n)_{h_1\v,h_2\v})_{h_3\v, h_4\v}.$$
Hence, the sum in (\ref{eqn5}) corresponding to
$f\in F$ is the negative of the sum in (\ref{eqn5})
corresponding to $\mu(f)$.
This completes the proof of Theorem \ref{main}.

\section*{Appendix: Strong homotopy Lie bialgebras}

In this appendix, we give the definition of
strong homotopy Lie (super) bialgebras.

A strong homotopy Lie bialgebra in the category
of super vector spaces over a field $\Bbbk$ of characteristic $0$
consists of the data of:
\begin{itemize}
\item
a super vector space $V= V_+ \oplus V_-$ over $\Bbbk$;
\item
an odd linear map $\a_{m,n}: (\Pi V)^{\otimes n} \to
(\Pi V)^{\otimes m}$ 
for each ordered pair of positive integers $m,n$.
\end{itemize}
These data are required to satisfy the following 
conditions:
\begin{itemize}
\item
$\a_{m,n}$ is commutative and cocommutative for all $m,n$;
\item
one has:
$$ \sum_{T\in T(m,n)} \a_{o(t(T)),i(t(T))} \circ_T 
\a_{o(s(T)),i(s(T))} = 0 $$
for all $m,n$.
\end{itemize}
We remind that the set $T(m,n)$ of flowcharts
was defined in Section 1, and for each $T\in T(m,n)$,
$$ \a_{o(t(T)),i(t(T))} \circ_T \a_{o(s(T)),i(s(T))} :
(\Pi V)^{\otimes n} \to (\Pi V)^{\otimes m} $$
is the composition of $\a_{o(t(T)),i(t(T))}$ and
$\a_{o(s(T)),i(s(T))}$ according to $T$.

The above definition comes from the theory
of Koszul duality of dioperads, cf. \cite{G}. The reader 
may refer to \cite{M} for further details.

\section*{Acknowledgments}
The author was partially supported by the NSF grant  DMS-0401509.

\end{document}